\documentclass[11pt]{amsart}

\usepackage{amsmath,amssymb,latexsym,amsfonts,amscd}

\theoremstyle{plain}

\newtheorem{theorem}{Theorem}[section]

\newtheorem{proposition}[theorem]{Proposition}

\newtheorem{lemma}[theorem]{Lemma}

\newtheorem{corollary}[theorem]{Corollary}

\theoremstyle{definition}

\theoremstyle{remark}

\newtheorem{remark}[theorem]{Remark}

\newcommand{\Homr}{\op{Hom}^{\mathrm{r}}}

\newcommand{\F}[1][]{F^{#1*}}                      

\newcommand{\D}[1][]{D^{(#1)}}                     

\newcommand{\op}[1]{\operatorname{#1}}

\renewcommand{\bar}[1]{\overline{#1}}

\renewcommand{\phi}{\varphi}

\renewcommand{\theta}{\vartheta}

\renewcommand{\epsilon}{\varepsilon}

\newcommand{\Ann}{\op{Ann}}         

\newcommand{\End}{\op{End}}

\newcommand{\id}{\op{id}}

\newcommand{\tensor}{\otimes}         

\renewcommand{\to}[1][]{\xrightarrow{\ #1\ }}

\newcommand{\usc}[1][m]{\underline{\phantom{#1}}}

\newcommand{\defeq}{\stackrel{\scriptscriptstyle \op{def}}{=}}

\begin{document}

\title[Generators of $D$--modules]{Generators of $D$--modules
in positive characteristic}
\author[J. Alvarez-Montaner]{Josep Alvarez-Montaner}
\thanks{Research of the first author partially supported by a
Fulbright grant and the Secretar\'{\i}a de Estado de Educaci\'on y
Universidades of Spain and the European Social Funding}
\address{Departament de Matem\`atica Aplicada I\\
Universitat Polit\`ecnica de Catalunya\\ Avinguda Diagonal 647,
Barcelona 08028, Spain} \email{Josep.Alvarez@upc.es}

\author[M. Blickle]{Manuel Blickle}
\thanks{The second author received funding from the DFG Schwerpunkt
``Globale Methoden in der Komplexen Geometrie''}
\address{FB6 Mathematik, Universit\"at Duisburg--Essen, 45117 Essen,
Germany}
\email{manuel.blickle@uni-essen.de}

\author[G. Lyubeznik]{Gennady Lyubeznik}
\thanks{The third author greatfully acknowledges NSF support}
\address{Department of Mathematics, University of Minnesota, Minneapolis,
MN 55455, USA} \email{gennady@math.umn.edu}

\begin{abstract}

Let $R=k[x_1, \dots ,x_d]$ or $R=k[[x_1,\dots,x_d]]$ be either a
polynomial or a formal power series ring in a finite number of
variables over a field $k$ of characteristic $p>0$ and let $D_{R|k}$
be the ring of $k$-linear differential operators of $R$. In this
paper we prove that if $f$ is a non-zero element of $R$ then
$R_f$, obtained from $R$ by inverting $f$, is generated as a
$D_{R|k}$--module by $ \frac{1}{f}$. This is an amazing fact
considering that the corresponding characteristic zero statement
is very false. In fact we prove an analog of this result for a
considerably wider class of rings $R$ and a considerably wider
class of $D_{R|k}$-modules.

\end{abstract}

\maketitle

\section{Introduction}
Let $k$ be a field and let $R=k[x_1,\dots, x_d]$, or
$R=k[[x_1,\dots, x_d]]$ be either a ring of polynomials or formal
power series in a finite number of variables over $k$. Let $D_{R|k}$
be the ring of $k$-linear differential operators on $R$. For every
$f\in R$, the natural action of $D_{R|k}$ on $R$ extends uniquely to
an action on the localization $R_f$ via the standard quotient
rule. Hence $R_f$ acquires a natural structure of
$D_{R|k}$--module. It is a remarkable fact that $R_f$ has finite
length in the category of $D_{R|k}$-modules. This fact has been
proven in characteristic $0$ by Bernstein \cite[Cor. 1.4]{Be72} in
the polynomial case and by Bj\"ork \cite[Thm. 2.7.12, 3.3.2]{Bj79}
in the formal power series case. In positive characteristic the
polynomial case was established by B\o gvad \cite[Prop.
3.2]{Bog.DmodBorel} and the formal power series case by Lyubeznik
\cite[Thm. 5.9]{Lyub}. Consequently the ascending chain of
$D_{R|k}$-submodules $$D_{R|k} \cdot \frac{1}{f} \subseteq D_{R|k} \cdot
\frac{1}{f^2} \subseteq \dots \subseteq D_{R|k} \cdot \frac{1}{f^s}
\subseteq \dots \subseteq R_f$$ stabilizes, i.e. $R_f$ is
generated by $\frac{1}{f^i}$ for some $i$. This paper is motivated
by the natural question: {\it What is the smallest $i$ such that
$\frac{1}{f^i}$ generates $R_f$ as a $D_{R|k}$--module?}

If $k$ is a field of characteristic zero and $f\in R$ is a
non-zero polynomial, the answer to this question is known:
Theorem 1' of \cite{Be72} shows that there exists of a monic
polynomial $b_f(s)\in k[s]$ and a differential operator $Q(s) \in
D_{R|k}[s]$ such that $$Q(s)\cdot f^{s+1}=b_f(s)\cdot f^s$$ for every
$s$. The polynomial $b_f(s)$ is called the Bernstein-Sato
polynomial of $f$. Let $-i$ be the negative integer root of
$b_f(s)$ of greatest absolute value (it exists since $-1$ is
always a root of $b_f(s)$). Then, $b_f(s)\ne 0$ for any integer
$s<-i$, hence $f^s\in D_{R|k}\cdot f^{s+1}$ implying $\frac{1}{f^s}\in
D_{R|k}\cdot\frac{1}{f^i}$ for all $s>i$. In particular, $R_f$ is
$D_{R|k}$-generated by $\frac{1}{f^i}$ and, as is shown in
\cite[Lem. 1.3]{Wa03}, it cannot be generated by $\frac{1}{f^j}$
for $j<i$. This gives a complete answer to our question in
characteristic zero.

For example, consider the polynomial $f=x_1^2+\dots+x_{2n}^2$.
Then, the functional equation
$$\frac{1}{4}(\frac{\partial^2}{\partial x_1^2} +
\dots + \frac{\partial^2}{\partial x_{2n}^2} )\cdot f^{s+1} =
(s+1)(s+n) \cdot f^s$$ shows that the Bernstein-Sato polynomial is
$b_f(s)=(s+1)(s+n)$ \cite[Cor. 3.17]{Yano}. Hence $R_f$ is
$D_{R|k}$-generated by $ \frac{1}{f^n}$ but it cannot be generated
by $ \frac{1}{f^j}$ for $j<n$.


The goal of this paper is to prove the following amazing result.

\begin{theorem}\label{T1}
Let $R=k[x_1,\dots,x_d]$ or $R=k[[x_1,\dots,x_d]]$ be either a
ring of polynomials or formal power series over a field $k$ of
characteristic $p>0$ and let $f\in R$ be non-zero. Then $R_f$ is
$D_{R|k}$-generated by $\frac{1}{f}$.
\end{theorem}

In fact we prove this result for a considerably wider class of rings $R$
and a considerably wider class of $D_{R|k}$-modules.

In Section 2 we have collected for the reader's convenience some
basic (and not so basic) facts about $D_{R|k}$-modules in
characteristic $p>0$. These are needed in Section 4 and 5, whereas
Section 3 can be read with minimal prior exposure to $D_{R|k}$--module
theory.

In Section \ref{seq} we introduce a chain of ideals associated to
an element $f$ of a regular $F$-finite ring $R$ of characteristic
$p>0$. These ideals are of considerable interest in their own
right and are likely to become useful in other contexts as well,
especially in the theory of tight closure. The crucial fact is
that this chain of ideals stabilizes if and only if $R_f$ is
$D_{R|k}$-generated by $\frac{1}{f}$ (Corollary \ref{coro}). This
yields a very elementary proof that $R_f$ is $D_{R|k}$-generated by
$\frac{1}{f}$ in the special case that $R$ is a polynomial ring
over a field (Theorem \ref{poly} and Corollary \ref{polycoro}).

In Section 4 we extend our results to all regular rings $R$ that
are of finite type over an $F$-finite regular local ring of
characteristic $p>0$ and to all finitely generated unit
$R[F]$-modules (Theorem \ref{thm.Main1} and Corollary
\ref{coro.main}). The proof uses considerably more advanced tools
than the elementary proof of Section 3 for the polynomial ring.
Namely, Frobenius descent and Lyubeznik's theorem \cite[Cor.
5.8]{Lyub} to the effect that $R_f$ has finite length in the
category of $D_{R|k}$-modules are used. This result implies that the
chain of ideals constructed in Section 3 stabilizes for this class
of rings. These ideals squarely belong to commutative algebra and
it is quite remarkable that the only available proof that they
stabilize requires the use of $D_{R|k}$-modules!

In Section 5 we extend our results to the case of regular algebras
$R$ of finite type over a formal power series ring
$k[[x_1,\dots,x_d]]$ where $k$ is an arbitrary field of
characteristic $p>0$ (the case that $k$ is perfect
is covered by the results of Section
4).

This paper combines and generalizes the results of preprints
\cite{MontLyub} and \cite{Bli.prep}.

\section{Rings of differential operators and modules over them
in characteristic $p>0$}\label{App} In this purely expository
section we have collected some basic
facts which are needed in the following sections. Throughout this
section $R$ is a commutative ring containing a field of
characteristic $p>0$.

\subsection{Definition and elementary properties}
The differential operators $\delta:R\to R$ of order $\le n$, where
$n$ is a non-negative integer, are defined inductively as follows
(cf. \cite[\S 16.8]{GD67}). A differential operator of order 0 is
just the multiplication by an element of $R$. A differential operator
of order $\le n$ is an additive map $\delta:R\to R$ such that for
every $r\in R$, the commutator $[\delta,\tilde r]=\delta\circ
\tilde r-\tilde r\circ \delta$ is a differential operator of order
$\le n-1$ where $\tilde r:R\to R$ is the multiplication by $r$.
The sum and the composition of two differential operators are
differential operators, hence the differential operators form a
ring which is a subring of End$_{\mathbb Z}R$. We denote this ring
$D_{R}$.

If $k\subseteq R$ is a subring, we denote by $D_{R|k}\subseteq
D_R$ the subring of $D_R$ consisting of all those differential
operators that are $k$-linear. Since every additive map $R\to R$
is $\mathbb Z/p\mathbb Z$-linear, $D_R=D_{R|\mathbb Z/p\mathbb
Z}$, i.e. $D_R$ is a special case of $D_{R|k}$.  The ring
homomorphism $R\to D_{R|k}$ that sends $r\in R$ to the
multiplication by $r$ makes $R$ a subring of $D_{R|k}$.

By a $D_{R|k}$--module we always mean a left $D_{R|k}$--module.
For example, $R$ with its natural $D_{R|k}$-action is a
$D_{R|k}$--module. If $M$ is a $D_{R|k}$--module and $S\subset R$
is a multiplicatively closed set, then $M_S$ has a unique
$D_{R|k}$--module structure such that the natural localization map
$M\to M_S$ is a $D_{R|k}$--module homomorphism \cite[Ex.
5.1]{Lyub}. In particular, $R_f$ carries a natural
$D_{R|k}$--module structure for every $f\in R$.

Every differential operator $\delta\in D_{R}$ of order $\leq
p^s-1$ is $R^{p^s}$-linear, where $R^{p^s}\subseteq R$ is the
subring consisting of all the $p^s$-th powers of all the elements
of $R$ \cite[1.4.8a]{Yeku}. In other words, $D_R$ is a subring of
the ring $$ \bigcup_{s}D^{(s)}_R $$ where
$D^{(s)}_R=\text{End}_{R^{p^s}}(R)$. In particular, this implies
that if $k$ is not just $\mathbb Z/p\mathbb Z$, but any perfect
subfield of $R$, i.e. $k\subseteq R^{p^s}$ for every $s$, then
every $\delta\in D_R$ is $k$-linear, i.e. $D_R=D_{R|k}$.

Let $k[R^{p^s}]$ be the $k$-subalgebra of $R$ generated by the
$p^s$-th powers of all the elements of $R$. If $R$ is a finite
$k[R^p]$--module, then \cite[1.4.9]{Yeku}
\begin{equation}
D_{R|k}=\bigcup_{s}\text{End}_{k[R^{p^s}]}(R).
\end{equation}

The ring $R$ is called {\it $F$-finite} if $R$ is a finitely
generated $R^p$--module. Since $D_R=D_{R|\mathbb Z/p\mathbb Z}$
and $\mathbb Z/p\mathbb Z[R^{p^s}]=R^{p^s}$, it  follows that if
$R$ is $F$-finite, then
\begin{equation}
  D_R=\bigcup_{s}D^{(s)}_R.
\end{equation}

If $R=k[x_1,\dots,x_d]$ or $R=k[[x_1,\dots,x_d]]$ is either the
ring of polynomials or the ring of formal power series over $k$, then
$D_{R|k}$ is the ring extension of $R$ generated by the
differential operators
$D_{t,i}=\frac{1}{t!}\frac{\partial^{t}}{\partial x_i^{t}}$ where
$\frac{\partial^{t}}{\partial x_i^{t}}$ is the $t$-th $k$-linear
partial differentiation with respect to $x_i$, i.e.
$D_{t,i}(x_i^s)=0$ if $s<t$ and
$D_{t,i}(x_i^s)=\binom{t}{s}x_i^{s-t}$ if $s\geq t$ \cite[\S
16.11]{GD67}, \cite[Ex. 5.3d,e]{Lyub}. If $k$ is perfect, i.e.
$\bigcup_{s}D^{(s)}_R=D_R=D_{R|k}$, then $D_R^{(s)}$ is the ring
extension of $R$ generated by the operators $D_{t,i}$ with $t<p^s$.

\subsection{Frobenius descent} The exposition in this subsection is based
on Chapter 3.2 of \cite{Bli.PhD}. We state and prove the basic result
but refer to \cite{Bli.PhD} for all the straightforward (but
tedious) compatibilities one has to check.

Frobenius descent has been used by a number of authors; see for
example S.P. Smith \cite{SmithSP.diffop,SmithSP:DonLine}, B.
Haastert \cite{Haastert.DiffOp,Haastert.DirIm} and R. B{\o}gvad
\cite{Bog.DmodBorel}. Its precursor is Cartier descent\footnote{It
states that $F^*$ is an equivalence between the category of
$R$--modules and the category of modules with integrable
connection and $p$--curvature zero. The inverse functor of $\F$ on
a module with connection $(M,\nabla)$ is in this case given by
taking the horizontal sections $\ker \nabla$ of $M$. As an
$R$--module with integrable connection and $p$--curvature zero is
nothing but a $\D[1]_R$--module, Cartier descent is just the case
$e=1$ of Proposition \ref{prop.FrobDesc}.} as described, for example, by
N.~M.~Katz
\cite[Thm. 5.1]{Katz}. A big generalization has recently been given by
P.~Berthelot
\cite{Ber.FrobDesc}.

In the basic form used here Frobenius descent is based on the fact
that a ring $R$ is Morita equivalent to the algebra of $n\times n$
matrices with entries in $R$. That is $R$ and Mat$_{n\times n}(R)$
have equivalent module categories. Namely, Mat$_{1\times n}(R)$,
the rows of length $n$ with entries from $R$ (resp. Mat$_{n\times
1}(R)$, the columns of length $n$ with entries from $R$) is an
$R$-Mat$_{n\times n}(R)$-bimodule (resp. a Mat$_{n\times
n}(R)$-$R$ bimodule) and the maps
$${\rm Mat}_{1\times n}(R)\otimes_{\rm Mat_{n\times n}(R)}{\rm
Mat}_{n\times 1}(R)\to R$$
$${\rm Mat}_{n\times 1}(R)\otimes_R{\rm Mat}_{1\times
n}(R)\to {\rm Mat}_{n\times n}(R)$$ that send $A\otimes B$ to the
matrix product $AB$ are isomorphisms, hence the functors
$${\rm Mat}_{1\times n}(R)\otimes_{{\rm Mat}_{n\times n}(R)}(\usc):{\rm
Mat}_{n\times n}(R){\rm -mod}\to R{\rm -mod}$$
$${\rm Mat}_{n\times 1}(R)\otimes_R(\usc):R{\rm -mod}\to {\rm
Mat}_{n\times n}(R){\rm -mod}$$ are inverses of each other and
establish an equivalence of categories.

Let $R^{(s)}$ be the abelian group of $R$ with the usual left
$D_R^{(s)}$--module structure (and hence the usual left
$R$--structure) and with the right $R$--module structure defined
by $r'r=r^{p^s}r'$ for all $r\in R$ and $r'\in R^{(s)}$. Thus
$R^{(s)}$ is a $D_R^{(s)}$--$R$--bimodule. We define a structure
of $R$--$D_R^{(s)}$--bimodule on $\Homr_R(R^{(s)},R)$ where
$\Homr$ denotes the homomorphisms in the category of right
$R$-modules as follows.  For $\delta \in D_R^{(s)}$, $\phi \in
\Homr_R(R^{(s)},R)$ and $r \in R$ the product $r \cdot \phi \cdot
\delta$ is the composition $\tilde r \circ \phi \circ \delta$
where $\delta$ acts on the left on $R^{(s)}$ and $\tilde r$ is the
multiplication by $r$ on $R$. The identification of $D^{(s)}_R$
with $\Homr_R(R^{(s)},R^{(s)})$ shows that this composition
$\tilde r \circ \phi \circ \delta$ is indeed in
$\Homr_R(R^{(s)},R)$. Thus we have functors
\begin{gather*}
F^{s*}(\usc)\defeq R^{(s)}\otimes_R\usc: \quad R{\rm -mod}\to
D_R^{(s)}{\rm-mod} \\
T^{s}(\usc) \defeq \Homr_R(R^{(s)},R) \tensor_{\D[s]_R}\usc: \quad
D^{(s)}_R{\rm -mod}\to R{\rm -mod}
\end{gather*}
the first of which is called the ($s$-fold) Frobenius functor on
$R$--modules.

\begin{proposition}[Frobenius Descent]\label{prop.FrobDesc}
If\/ $R$ is regular and\/ $F$--finite, the functors $F^{s*}(\usc)$
and $T^{s}(\usc)$ are inverses of each other. Consequently they
induce  an equivalence between the category of\/ $R$--modules and
the category of $\D[s]_R$--modules.
\end{proposition}

\begin{proof}
To show that $F^{s*}(\usc)$ and $T^{s}(\usc)$ are inverses of each
other it is enough to show that the natural map
\[
\Phi:R^{(s)} \tensor_R \Homr_R(R^{(s)},R) \to \D[s]_R
\]
given by sending $a \tensor \phi$ to the composition
\[
R^{(s)} \to[\phi] R \to[\tilde a] R \to[\rm id] R^{(s)}
\]
(where id is the identity map on the underlying abelian group of $R$
and $R^{(s)})$ and the natural map
\[
  \Psi:\ \Homr_R(R^{(s)},R) \tensor_{\D[s]_R} R^{(s)} \to R
\]
given by sending $\phi \tensor a$ to $\phi(a)$ are both ring
isomorphisms.

They are isomorphisms if and only if they are isomorphisms
locally. Since $R$ is regular, $R^{(s)}$ is a locally free right
$R$--module of finite rank \cite{Kunz}. Once an $R$-basis of
$R^{(s)}$ is fixed, we may view $R^{(s)}$ as the set of coordinate
rows of the elements of $R^{(s)}$ with respect to this basis and
$\Homr(R^{(s)},R)$ as the set of the coordinate columns of the
elements of $\Homr(R^{(s)},R)$ with respect to the dual basis, and
$\D[s]_R$ is just the matrix algebra over $R$, so we are done by
Morita duality between $R$ and Mat$_{n\times n}$, as described above.
\end{proof}

\begin{remark}
For a more explicit description of $T^s$ let $J_s$ be the left
ideal of $\D[s]_R$ consisting of all $\delta$ such that
$\delta(1)=0$. Then it is shown in \cite[Prop. 3.12]{Bli.PhD} that
$T^s(M) \cong \Ann_M J_s \subseteq M$.
\end{remark}

Proposition \ref{prop.FrobDesc} implies that the categories of
$\D[s]_R$--modules for all $s$ are equivalent since each single
one of them is equivalent to $R$--mod. The functor giving the
equivalence between $D_R^{(t)}$--mod and $\D[t+s]_R$--mod is, of
course, $\F[s]$. Concretely, to understand the $\D[t+s]_R$--module
structure on $\F[s]M$ for some $\D[t]_R$--module $M$, we write $M
\cong \F[t]N$ for $N=T^t(M)$. Then
$\F[s]M=\F[(t+s)]N=R^{(t+s)}\tensor_R N$ carries obviously a
$\D[t+s]_R$--module structure with $\delta \in \D[t+s]_R$ acting
via $\delta \tensor \id_N$.

Since the union $\bigcup_s \D[s]_R$ is just the ring of
differential operators $D_R$ of $R$ this implies the following
proposition (after the obvious compatibilities are checked, which
is straightforward and carried out in \cite[Chapter
3.2]{Bli.PhD}):

\begin{proposition}\label{prop.FrobDesc4Dmod}
Let\/ $R$ be regular and\/ $F$--finite. Then\/ $\F[s]$ is an
equivalence of the category of\/ $D_R$--modules with itself.
\end{proposition}

\subsection{Unit R[F]-modules}
We denote by $R[F]$ the ring extension of $R$ generated by a
variable $F$ subject to relations $Fr=r^pF$ for all $r\in R$.
Clearly, a (left) $R[F]$--module is an $R$--module $M$ with a map
of additive abelian groups $F:M\to M$ such that $F\circ \tilde
r=\tilde r^p\circ F$ where $\tilde r:R\to R$ is the multiplication
by $r$. To every $R[F]$--module $M$ is associated the map of
$R$-modules $\theta_M:F^*M=R^{(1)}\otimes_RM\to M$ sending
$r\otimes m$ to $rF(m)$. The $R[F]$--module $M$ is called a
\emph{unit $R[F]$--module} if $\theta_M$ is bijective. Unit
$R[F]$-modules are called $F$-modules in \cite{Lyub}.

A unit $R[F]$--module $(M,\theta)$ carries a natural structure of
$\bigcup_{s}D^{(s)}_R$--module and hence also of $D_R$--module as $D_R$
is a subring of $\bigcup_{s}D^{(s)}_R=\bigcup_s \End_{R^{p^s}}(R)$.
Namely, set
$$\theta_s=F^{(s-1)*}(\theta_M^{-1})\circ F^{(s-2)*}(\theta_M^{-1})\circ
\dots\circ\theta_M^{-1}:M\to F^{s*}M.$$ Every $u\in
U_s=\text{End}_{R^{p^s}}(R)$ acts on $F^{s*}M=R^{(s)}\otimes_RM$
via $u\otimes_R\text{id}_{M}$. We let $u$ act on $M$ via
$\theta_s^{-1}\circ(u\otimes_R\text{id}_{M})\circ\theta_s$. This
action is well-defined, i.e. it depends only on $u$, but not on
the particular $s$ \cite[p. 116]{Lyub}.

\begin{lemma}\label{lem.ThetaDlinear}
  Let\/ $R$ be regular and\/ $F$--finite and let\/ $M$ be a
unit\/ $R[F]$--module. Then\/ $\theta_M: F^*M \to M$ is a map of\/
$D_R$--modules where the $D_R$-structure on $F^*M$ is due to
Theorem \ref{prop.FrobDesc}.
\end{lemma}

\begin{proof}
We omit a straightforward verification of this and instead refer
to \cite[Chapter 3.2]{Bli.PhD}
\end{proof}

The usual $D_{R}$--module structure on $R_f$ is induced, as above,
by the unit $R[F]$--module structure $F:R_f\to R_f$ sending $x\in
R_f$ to $x^p$ \cite[Ex. 5.2c]{Lyub}. The $R[F]$--module $R_f$ is
generated by $\frac{1}{f}$ because
$F^{s}(\frac{1}{f})=\frac{1}{f^s}$, i.e. $R_f$ is a finitely
generated unit $R[F]$--module (finitely generated unit
$R[F]$-modules are called $F$-finite modules in \cite{Lyub}).

\begin{theorem} \label{Lyubez}(\cite[Cor. 5.8]{Lyub})
Let $R$ be a regular finitely generated algebra over a commutative
Noetherian regular $F$-finite local ring $A$ of characteristic
$p>0$. A finitely generated unit $R[F]$--module $M$ has finite
length in the category of $D_R$-modules. In particular, $R_f$ with
its usual $D_R$--module structure has finite length in the
category of $D_R$-modules for every $f\in R$.
\end{theorem}

\section{A chain of ideals associated to an element of a
regular $F$-finite ring}\label{seq}

In this section $R$ is a regular and $F$--finite ring of
characteristic $p>0$. For a given element $f \in R$ we aim to
define a descending chain of ideals $I_s(f)$ indexed by the positive
integers.

For this let us first assume that $R$ is a free $R^{p^s}$--module.
Let $\{c^{p^s}_1,c^{p^s}_2,\dots\}\subset R^{p^s}$ be the string
of coordinates of $f\in R$ with respect to some $R^{p^s}$-basis of
$R$. We define $I_s(f)$ to be the ideal of $R$ generated by the
set $\{c_1,c_2,\dots\}$. This definition is independent of the
choice of basis because the coordinates of $f$ with respect to one
basis are linear combinations with coefficients from $R^{p^s}$ of
the coordinates of $f$ with respect to another basis, hence the
corresponding ideals are the same.

Since any regular $F$--finite ring $R$ is a finite locally free
$R^p$--module \cite{Kunz}, Spec $R$ is covered by a finite number
of open affines Spec $R_{r}$ such that $R_{r}$ is a free
$R^p_{r}$--module (and consequently $R_{r}$ is a free
$R_{r}^{p^s}$--module for every $s$). Hence we define $I_s(f)$ in
general by glueing the local ideals defined above. This is
possible due to the independence of the choice of basis in the
construction.

This section is devoted to the study of the ideals $I_s(f)$
leading to an elementary proof of our main result in the
polynomial case. But these ideals are interesting by themselves
and are likely to become important, for example in tight closure
theory.

A further consequence of $R$ being $F$-finite, is that the ring of
differential operators of $R$ is $\bigcup_s D_R^{(s)}$ where
$D_R^{(s)}={\rm End}_{R^{p^s}}(R)$, according to formula (2) of
Section 2. One has the following relationship between the ideals
$I_s(f)$ and differential operators.

\begin{lemma}\label{D_R^s(f)}
$D_R^{(s)}\cdot f=I_s(f)^{[p^s]}$ where $D_R^{(s)}\cdot
f\stackrel{\rm def}{=}\{\delta(f)|\delta\in D_R^{(s)}\}\subseteq
R$ and $I_s(f)^{[p^s]}$ is the ideal generated by the $p^s$-th
powers of the elements of $I_s(f)$, equivalently, by the $p^s$-th
powers of a set of generators of $I_s(f)$.
\end{lemma}

\begin{proof} Since $R$ is a finitely generated $R^{p^s}$--module,
$D_R^{(s)}={\rm End}_{R^{p^s}}R$ commutes with localization with
respect to any multiplicatively closed subset of $R^{p^s}$. Hence
we may assume that $R$ is a free $R^{p^s}$--module. In this case
$f=\sum_ic_i^{p^s}e_i$ where $\{e_1,e_2,\dots\}$ is an
$R^{p^s}$-basis of $R$ and $\{c_1^{p^s},c_2^{p^s},\dots\}$ are the
coordinates of $f$ with respect to this basis. Since
$\delta(f)=\sum_i c_i^{p^s}\delta(e_i)\in
I_s(f)^{[p^s]}=(c_1^{p^s},c_2^{p^s},\dots)$ for every $\delta\in
D_{R}^{(s)}$, we see that $D_{R}^{(s)}\cdot f\subseteq
I_s(f)^{[p^s]}$. Conversely, setting $\delta_i\in D_{R}^{(s)}$ to
be the $R^{p^s}$-linear map that sends $e_i$ to $1$ and $e_j$ to
$0$ for every $j\ne i$ we see that $\delta_i(f)=c_i^{p^s}$, i.e.
every generator of $I_s(f)^{[p^s]}$ is in $D_{R}^{(s)}$.
\end{proof}

\begin{lemma}\label{f^p}
$I_s(f)=I_{s+1}(f^p)$.
\end{lemma}

\begin{proof} It is enough to prove this after localization at every
maximal ideal of $R$, hence we can assume that $R$ is local in which case
$R$ is a free $R^p$-module. Since
$1\not\in \mathfrak mR$, where $\mathfrak m$ is the maximal
ideal of $R^p$,
Nakayama's lemma implies that we can take
$1$ to be part of a free
$R^p$-basis of
$R$, i.e. we may assume that
$R$ is a free
$R^{p}$--module on some basis $\{e_1,e_2,\ldots\}$ and
$e_1=1$.

Now let $\{\tilde e_j\}$ be an $R^{p^s}$-basis of $R$. Then the
set of all products $\{ e_{j,i}=\tilde e^p_je_i\}$ is an
$R^{p^{s+1}}$--basis of $R$. Raising the equality
$f=\sum_jc_j^{p^s}\tilde e_j$ to the $p$-th power we get
$f^p=\sum_j c_j^{p^{s+1}}\tilde e_j^p$. But $\tilde e_j^p=\tilde
e_j^p\cdot 1=\tilde e_j^p\cdot e_1=e_{j,1}$, hence the $(j,i)$-th
coordinate of $f^p$ with respect to the basis $\{e_{j,i}\}$ is
$c_j^{p^{s+1}}$ if $i=1$ and $0$ if $i\ne 1$. Hence $I_s(f)$ and
$I_{s+1}(f^p)$ are generated by the same elements. \end{proof}

\begin{lemma}\label{subset}
$I_s(f\tilde f)\subseteq I_s(f)I_s(\tilde f)\subseteq I_s(f)$ for
every $\tilde f\in R$.
\end{lemma}

\begin{proof} As before we may assume that $R$ is a free
$R^{p^s}$--module on some basis $\{e_1,e_2,\ldots\}$. Multiplying
the equalities $f=\sum_ic_i^{p^s}e_i$ and $\tilde
f=\sum_i\tilde c_i^{p^s}e_i$ we get $f\tilde f=\sum_{i,j}c_i^{p^s}\tilde
c_j^{p^s}e_ie_j$. Writing $e_ie_j=\sum_q\bar c_q^{p^s}e_q$,
substituting this into the preceding equality and collecting
similar terms we see that all the coordinates of $f\tilde f$ with
respect to the basis $\{e_i\}$ are linear combinations (with
coefficients from $R^{p^s}$) of the products $c_i^{p^s}\tilde
c_j^{p^s}$. This implies that $I_s(f\tilde f)$ is generated by
linear combinations (with coefficients from $R$) of the products
$c_ic_j$, hence $I_s(f\tilde f)\subseteq I_s(f)I_s(\tilde
f)\subseteq I_s(f)$. \end{proof}

\begin{lemma}\label{p^s-1}
$I_{s+1}(f^{p^{s+1}-1})\subseteq I_s(f^{p^s-1})$.
\end{lemma}

\begin{proof} $f^{p^{s+1}-1}=f^{p^{s+1}-p}f^{p-1}$ so
$I_{s+1}(f^{p^{s+1}-1})\subseteq I_{s+1}(f^{p^{s+1}-p})$ by Lemma
\ref{subset}. Since $f^{p^{s+1}-p}=(f^{p^s-1})^p$, we are done by
Lemma \ref{f^p}.\end{proof}

\begin{proposition}\label{stabilization}
The descending chain of ideals \[ I_1(f^{p-1})\supseteq
I_2(f^{p^2-1})\supseteq \dots\] stabilizes at $s$, i.e.
$I_{s}(f^{p^{s}-1})=I_{s+1}(f^{p^{s+1}-1})=I_{s+2}(f^{p^{s+2}-1})=\dots$,
if and only if there is $\delta\in D^{(s+1)}_R$ such that $\delta(
\frac{1}{f})=\frac{1}{f^p}.$
\end{proposition}

\begin{proof} Assume $I_{s}(f^{p^{s}-1})=I_{s+1}(f^{p^{s+1}-1})$.
On the other hand Lemma \ref{f^p} implies that
$I_{s}(f^{p^{s}-1})=I_{s+1}(f^{p^{s+1}-p})$. Hence
$$I_{s+1}(f^{p^{s+1}-p})=I_{s+1}(f^{p^{s+1}-1})$$ and consequently
$I_{s+1}(f^{p^{s+1}-p})^{[p^{s+1}]}=I_{s+1}(f^{p^{s+1}-1})^{[p^{s+1}]}$.
This implies that $D_R^{(s+1)}\cdot f^{p^{s+1}-p} =
D_R^{(s+1)}\cdot f^{p^{s+1}-1}$ by Lemma \ref{D_R^s(f)}. But
$f^{p^{s+1}-p}=\delta'(f^{p^{s+1}-p})$ where $\delta'=1\in
D_R^{(s+1)}$ so we have $f^{p^{s+1}-p}\in D^{(s+1)}_R\cdot
f^{p^{s+1}-p}$. Hence $$f^{p^{s+1}-p}\in D_R^{(s+1)}\cdot
f^{p^{s+1}-1}$$ i.e. there is $\delta\in D_R^{(s+1)}$ such that
$\delta(f^{p^{s+1}-1})=f^{p^{s+1}-p}$. Dividing this equality by
$f^{p^{s+1}}$ and considering that every $\delta\in D_R^{(s+1)}$
commutes with every element of $R^{p^{s+1}}$ we get $\delta(
\frac{f^{p^{s+1}-1}}{f^{p^{s+1}}})=\frac{f^{p^{s+1}-p}}{f^{p^{s+1}}}$,
i.e. $\delta (\frac{1}{f})=\frac{1}{f^p}.$

Conversely, assume there is $\delta\in D^{(s+1)}_R$ such that
$\delta(\frac{1}{f})=\frac{1}{f^p}.$ Multiplying this equality by
$f^{p^{s+1}}$ we get $\delta(f^{p^{s+1}-1})=f^{p^{s+1}-p}$. This
implies that $$\hskip -2.5cm D^{(s+1)}_R \cdot
f^{p^{s+1}-p}=D^{(s+1)}_R \cdot(\delta( f^{p^{s+1}-1}))=$$
$$\hskip 3.1cm =(D^{(s+1)}_R\cdot \delta) (f^{p^{s+1}-1})\subseteq
D^{(s+1)}_R\cdot f^{p^{s+1}-1}.$$

\noindent Hence $I_{s+1}(f^{p^{s+1}-p})^{[p^{s+1}]}\subseteq
I_{s+1}(f^{p^{s+1}-1})^{[p^{s+1}]}$ by Lemma \ref{D_R^s(f)}. As is shown
in the paragraph after next, this implies
$I_{s+1}(f^{p^{s+1}-p})\subseteq I_{s+1}(f^{p^{s+1}-1})$. Now Lemma
\ref{f^p} implies
$I_s(f^{{p^s}-1})\subseteq I_{s+1}(f^{p^{s+1}-1})$ since
$f^{p^{s+1}-p}=(f^{p^s-1})^p$. This containment and
Lemma \ref{p^s-1} imply that
$I_s(f^{{p^s}-1})=I_{s+1}(f^{p^{s+1}-1})$.

We have proven that the existence of $\delta\in D^{(s+1)}_R$ such
that $\delta(\frac{1}{f})=\frac{1}{f^p}$ is equivalent to equality
$I_s(f^{{p^s}-1})=I_{s+1}(f^{p^{s+1}-1})$. It now follows that
$I_s(f^{{p^s}-1})=I_{s+1}(f^{p^{s+1}-1})$ is equivalent to
$I_s(f^{{p^s}-1})=I_{s'}(f^{p^{s'}-1})$ for all $s'>s$ because
every $\delta\in D^{(s+1)}_R$ automatically belongs to
$D^{(s')}_R$ for all $s'>s$.

This completes the proof of the lemma modulo the fact, used
in the paragraph before the preceding one, that if $\mathcal
I$ and $\mathcal J$ are two ideals of $R$ such that $\mathcal
I^{[p^{s+1}]}\subset \mathcal J^{[p^{s+1}]}$, then $\mathcal I\subset
\mathcal J$. We are now going to prove this fact. It is enough to prove
this locally, hence we can assume like in the proof of Lemma 3.2 that
$R$ is a free $R^{p^{s+1}}$-module and $e_1=1$ is one of
the free generators, i.e. $R=R^{p^{s+1}}\oplus M$ where $M$ is a free
$R^{p^{s+1}}$-module. Let $\varphi:R\to R^{p^{s+1}}$ be the map sending
$r\in R$ to $r^{p^{s+1}}\in R^{p^{s+1}}$. Clearly,
$\mathcal I^{[p^{s+1}]}=\varphi(\mathcal I)R$, hence $\mathcal
I^{[p^{s+1}]}\cap R^{p^{s+1}}=\varphi(\mathcal I)R\cap
R^{p^{s+1}}=(\varphi(\mathcal
I)\oplus\varphi(\mathcal I)M)\cap R^{p^{s+1}}=\varphi(\mathcal
I)$, and the same holds with
$\mathcal I$ replaced by $\mathcal J$. Hence upon taking the intersection
with $R^{p^{s+1}}$ the containment $\mathcal
I^{[p^{s+1}]}\subset \mathcal J^{[p^{s+1}]}$ implies $\varphi(\mathcal
I)\subset \varphi(\mathcal J)$ which implies $\mathcal I\subset \mathcal
J$ since $\varphi$ is a ring isomorphism.
\end{proof}

We note that the proof shows that the descending chain of ideals
stabilizes at the first integer $s$ such that
$I_s(f^{{p^s}-1})=I_{s+1}(f^{p^{s+1}-1})$.

\begin{corollary}\label{coro}
The chain of ideals $I_1(f^{p-1})\supseteq I_2(f^{p^2-1})\supseteq
\dots$ stabilizes if and only if $\frac{1}{f}$ generates $R_f$ as
a $D_R$--module.
\end{corollary}

\begin{proof} If $\frac{1}{f}$ generates $R_f$ as a $D_R$--module,
then $\frac{1}{f^p}\in D_R\cdot \frac{1}{f}$, i.e. there is
$\delta\in D_R$ such that $\frac{1}{f^p}=\delta(\frac{1}{f})$.
Since $\delta\in D^{(s+1)}_R$ for some $s$, the preceding
proposition shows that the chain of ideals stabilizes at $s$.

Conversely, if the chain of ideals stabilizes at $s$, the
preceding proposition implies that there exists $\delta\in
D^{(s+1)}_R$ such that $\delta (\frac{1}{f})=\frac{1}{f^p}.$ As is
proven in the course of the proof of Lemma \ref{f^p}, locally the
element $1$ can always be taken as one of the elements of an
$R^{p}$-basis of $R$, hence $R/R^{p}$ is a finite locally free,
hence projective, $R^{p}$--module. Thus the natural surjection
$R\to R/R^{p}$ splits, i.e. there exists an $R^{p}$--module
isomorphism $R\cong R^{p}\oplus R/R^{p}$. Let $\delta'\in
D_R^{(s+2)}$ be defined by $\delta'(x^p\oplus y)=\delta(x)^p$ for
all $x\in R$ (i.e. $x^p\in R^p$) and $y\in R/R^p\subseteq R$. Then
$\delta'(\frac{1}{f^p})=(\delta(\frac{1}{f}))^p=(\frac{1}{f^p})^p=
\frac{1}{f^{p^2}}$, i.e. $\frac{1}{f^{p^2}}\in
D_R\cdot\frac{1}{f}$. Thus we have shown for any $f$ that
$\frac{1}{f^p}\in D_R\cdot\frac{1}{f}$ implies that
$\frac{1}{f^{p^2}}\in D_R\cdot\frac{1}{f}$. Hence
$\frac{1}{f^{p^s}}\in D_R\cdot\frac{1}{f}$ for every $s$, by
induction on $s$. But the set $\{\frac{1}{f^{p^s}}\}$ generates
$R_f$ as an $R$--module.\end{proof}

\emph{Open Question.} Let $R$ be a regular $F$-finite ring of
characteristic $p>0$ and let $f\in R$ be an element. Does the
chain of ideals $I_1(f^{p-1})\supseteq
I_2(f^{p^2-1})\supseteq\dots$ stabilize? Equivalently, is $R_f$
generated by $\frac{1}{f}$ as a $D_R$--module?

\medskip

For an arbitrary $F$-finite ring $R$ the question is open. But in the
next section we show that for a large class of regular $F$-finite rings
the answer is \emph{yes}!
First however, we give an elementary treatment for $R$ a polynomial ring.

\subsection{The case of a polynomial ring}
We will use multi-index notation in the case of the polynomial ring
$R=k[x_1,
\dots ,x_n]$, where $k$ is a field of characteristic $p>0$. A
differential operator $\delta\in D_R$ will be written in the right
normal form, i.e. $\delta =\sum a_{\alpha \beta} \hskip 1mm
{{x}^\alpha} \hskip 1mm D_\beta,$ where ${{x}^{\alpha}}$ will
stand for the monomial ${{x}^{\alpha}}:= x_1^{\alpha_1}\cdots
x_n^{\alpha_n}$, ${D_{\beta}}$ will denote the differential
operator ${D_{\beta}}:= D_{\beta_1,1}\cdots D_{\beta_n,n}$ and all
but finitely many $a_{\alpha \beta}\in k$ are zero.

\begin{theorem}\label{poly}
Let $R=k[x_1,\dots,x_d]$ be a polynomial ring in $x_1,\dots, x_d$
over a perfect field $k$ of characteristic $p>0$ and let $f\in R$
be any element. The chain of ideals $I_1(f^{p-1})\supseteq
I_2(f^{p^2-1})\supseteq\dots$ stabilizes.
\end{theorem}

\begin{proof} The monomials $\{ x^{\alpha}= x_1^{\alpha_1}\cdots
x^{\alpha_d}_d \hskip 1mm | \hskip 1mm 0\leq \alpha_i\leq p^s-1\}$
form an $R^{p^s}$-basis of $R$. Let
$f^{p^s-1}=\sum_{\alpha}c_{\alpha}^{p^s} x^{\alpha}$, so that
$I_s(f^{p^s-1})$ is generated by the set $\{c_{\alpha}\}$. No
monomials on the right side of the equation
$f^{p^s-1}=\sum_{\alpha}c_{\alpha}^{p^s} x^{\alpha}$ get cancelled
as a result of reducing similar terms, hence deg $f^{p^s-1}\geq
{\rm deg}\hskip 1mm c_{\alpha}^{p^s}$ for every $ \alpha$. This
inequality translates to $(p^s-1)\hskip 1mm{\rm deg}\hskip 1mm f
\geq p^s \hskip 1mm {\rm deg} \hskip 1mm c_{\alpha}$, i.e. ${\rm
deg}\hskip 1mm c_{\alpha}\leq \frac{p^s-1}{p^s}\hskip 1mm{\rm deg}
\hskip 1mm f$. Hence deg $c_{\alpha}<{\rm deg}\hskip 1mm f$. Thus
the ideals $I_s(f^{p^s-1})$ for every $s$ are generated by
polynomials of degrees less than deg $f$ which is independent of
$s$. The set of polynomials of degrees less than deg $f$ is a
finite dimensional $k$-vector space and the intersections of the
ideals $I_s(f^{p^s-1})$ with this vector space form a descending
chain of subspaces which stabilizes because the space is finite
dimensional. Hence $I_1(f^{p-1})\supseteq
I_2(f^{p^2-1})\supseteq\dots$ stabilizes. \end{proof}

\begin{corollary}\label{polycoro}
Let $R=k[x_1,\dots,x_d]$ be a polynomial ring in $x_1,\dots, x_d$
over an arbitrary field $k$ of characteristic $p>0$ and let $f\in R$
be any element. $R_f$ is
generated by $\frac{1}{f}$ as a $D_{R|k}$--module.
\end{corollary}

\begin{proof} If $k$ is perfect, we are done by Corollary \ref{coro} and
Theorem \ref{poly}. In the general case, let $K$
be the perfect closure of $k$. Since $K$ is perfect, there is a
differential operator $\delta=\sum a_{\alpha \beta} {{x}^\alpha} D_\beta$
with coefficients $a_{\alpha \beta}\in K$ such that
$\delta(\frac{1}{f})= \frac{1}{f^p}$. This is equivalent to the
fact that a system of finitely many linear equations with
coefficients in $k$ has solutions in $K$, where the non-zero
coefficients $a_{\alpha\beta}$ of $\delta$ are thought of as the
unknowns of the system. (For example, if $f=x_1$, we may be looking
for a solution in the form $\delta=a D_{p-1,1}$, so we get an
equation $\delta(\frac{1}{x})=\frac{1}{x^p}$. Since
$\delta(\frac{1}{x_1})=a \frac{1}{x_1^p}$, the corresponding
linear system is just one equation $a=1$.) The system has a
solution in $K$, namely, the coefficients of $\delta$. Hence it is
consistent, so it must have a solution in $k$ because the
coefficients of the linear system are in $k$ (they depend only on the
coefficients of $f$). So there is a differential operator $\delta'$ with
coefficients in
$k$ such that
$\delta'(\frac{1}{f})= \frac{1}{f^p}$.
\end{proof}

We conclude this section with an example. Let
$R=k[x_1,x_2,x_3,x_4]$ where $k$ is a field of characteristic
$p>0$ and let $f=x_1^2+x_2^2+x_3^2+x_4^2 $. In characteristic $0$,
as is pointed out in the Introduction, $\frac{1}{f^2}$ does not
belong to the $D_R$-submodule of $R_f$ generated by $\frac{1}{f}$.
But in characteristic $p>0$ we are going to find a differential operator
$\delta\in D_R$ such that $\delta(\frac{1}{f})=\frac{1}{f^p}$ just by
investigating the monomials appearing in $f^{p-1}$.

\begin{itemize}
\item If $4$ divides $p-1$, then $f^{p-1}$ contains the term $a_\alpha
{x}^\alpha$ where
\[
\textstyle a_\alpha= \frac{(p-1)!}{(\frac{p-1}{4}!)^4}\ne 0 \text{
and }
\alpha=(\frac{p-1}{2},\frac{p-1}{2},\frac{p-1}{2},\frac{p-1}{2}).
\]
\item If $4$ does not divide $p-1$, then $f^{p-1}$ contains the term
$a_\alpha
{x}^\alpha$ where
\[
  \textstyle a_\alpha=
\frac{(p-1)!}{(\frac{p+1}{4}!)^2(\frac{p-3}{4}!)^2}\ne 0\text{ and }
  \alpha=(\frac{p+1}{2},\frac{p+1}{2},\frac{p-3}{2},\frac{p-3}{2}).
\]
\end{itemize}

Notice that $\frac{1}{a_\alpha} \hskip 1mm D_\alpha(f^{p-1})=1$
because all other monomials appearing in $f^{p-1}$ contain some
$x_i$ raised to a power smaller than the power of
$\frac{\partial_i}{\partial x_i}$ in $D_\alpha$ hence $D_\alpha$
annihilates all other monomials. The differential operator
$\delta= \frac{1}{a_\alpha} \hskip 1mm D_\alpha$ commutes with
$f^p$ so, dividing the equation $\delta(f^{p-1})=1$ by $f^p$ we
get the desired result.

\section{The case of a regular finitely generated algebra
over an $F$-finite regular local ring}\label{mainsec}

Here we prove the central result of our paper using the techniques
surveyed in Section 2.

\begin{theorem}\label{thm.Main1}
Let $R$ be a regular finitely generated algebra over an $F$-finite
regular local ring of characteristic $p>0$. Let $f \in R$ be a
nonzero element. Then the $D_R$--module $R_f$ is generated by
$\frac{1}{f}$.
\end{theorem}

\begin{proof}[Proof]
For any $D_R$--submodule $M \subseteq R_f$ we identify $F^*M$ with
its isomorphic image in $R_f$ via the natural $D_R$--module
isomorphism $\theta: F^*R_f \to R_f$ of Lemma \ref{lem.ThetaDlinear}
(with $R_f$ viewed as a unit $R[F]$-module given by the map $F:R_f\to
R_f$ that sends every $x\in R_f$ to $x^p$). Then
$F^*M$ is
$R$-generated by the elements $m^p$ for $m \in M \subseteq R_f$. By
Frobenius descent (Proposition \ref{prop.FrobDesc}), $F^*M$ is a
$D_R$--submodule of $R_f$.

Let $M=D_R \cdot\frac{1}{f}$. We claim that $M \subseteq F^*M$.
Because $F^*M$ is a $D_R$--submodule of $R_f$, it is enough to
show that $\frac{1}{f} \in F^*M$. But $\frac{1}{f}\in M$ implies
$(\frac{1}{f})^p=\frac{1}{f^p}\in F^*M$, hence $f^{p-1}\cdot
\frac{1}{f^p}=\frac{1}{f}\in F^*M$. This proves the claim.

Now we get an ascending chain of $D_R$--submodules of $R_f$:
\begin{equation}\label{eq}
M \subseteq F^*M \subseteq F^{2*}M \subseteq F^{3*}M \subseteq
\ldots
\end{equation}

The fact that $\frac{1}{f} \in M$ implies
$\frac{1}{f^{p^s}}=F^s(\frac{1}{f})\in F^{s*}M$, hence the union
of the chain must be all of $R_f$. Thus it is enough to show that
$M=F^*M$ since then $M=F^{s*}M$ for all $s$, hence $M=R_f$ as
claimed. Assume otherwise, that is assume that the inclusion $M
\subsetneq F^*M$ is strict. Then \emph{all} the inclusions of
(\ref{eq}) must be strict since
$F^{s*}(\usc)=R^{(s)}\otimes_R(\usc)$ and $R^{(s)}$ is a
faithfully flat right $R$--module. But this contradicts the fact
that by Theorem \ref{Lyubez} the length of $R_f$ as a
$D_R$--module is finite.
\end{proof}

\begin{corollary}
Let $R$ be a regular finitely generated algebra over an $F$-finite
regular local ring of characteristic $p>0$. Let $f\in R$ be any
element. The descending chain of ideals $I_1(f^{p-1})\supseteq
I_2(f^{p^2-1})\supseteq\dots$ defined in the preceding section
stabilizes.
\end{corollary}

\begin{proof}
This follows from Corollary \ref{coro}.
\end{proof}

Theorem \ref{thm.Main1} also follows from the following more
general observation which was inspired by \cite[Prop.
15.3.4]{EmKis.Fcrys}, which in the notation of Theorem
\ref{thm.Main2} states that if $F^*M \subseteq M$ then $M$ is also
a unit $R[F]$--submodule.

\begin{theorem}\label{thm.Main2}
Let $R$ be a regular finitely generated algebra over an $F$-finite
regular local ring of characteristic $p>0$. Let $N$ be a finitely
generated unit $R[F]$--module. Suppose $M \subseteq N$ is a
$D_R$--submodule such that $M \subseteq F^*M$ (we identify $F^*M
\subseteq F^*N$ with  its image in $N$ via the structural
isomorphism $\theta: F^*N \to N$ of $N$). Then $M$ is a unit
$R[F]$--submodule.
\end{theorem}

\begin{proof}
$M$ being a unit $R[F]$--submodule of $N$ just means that the
inclusion $M \subseteq F^*M$ is in fact an equality. If the
inclusion is strict, then all the inclusions $F^{s*}M \subsetneq
F^{(s+1)*}M$ are strict as well because they are obtained by
tensoring $M \subseteq F^*M$ with the faithfully flat $R$--module
$R^{(s)}$. The resulting strictly increasing infinite chain
\[
  M \subsetneq F^*M \subsetneq F^{2*}M \subsetneq F^{3*}M
\subsetneq \cdots
\]
contradicts the finiteness of the length of $N$ as a
$D_R$--module.
\end{proof}

To obtain Theorem \ref{thm.Main1} from this just note that $M =
D_R \cdot \frac{1}{f}$ satisfies $M \subseteq F^*M$ and contains
the $R[F]$--module generator $\frac{1}{f}$ of $R_f$.

\medskip

An $R$--submodule $N_0$ of a unit $R[F]$--module $N$ is called a
\emph{root}, if $N_0$ is finitely generated as an $R$--module,
$N_0 \subseteq F^*N_0$ and $\bigcup_s F^{s*}N_0 = N$. The
existence of a root is equivalent to $N$ being finitely generated
as a unit $R[F]$--module \cite[Cor. 2.12]{Bli.PhD}.

\begin{corollary}\label{coro.main}
With the same assumptions as in Theorem \ref{thm.Main2}, if\/
$n_1,\ldots,n_t$ are generators of a root of a finitely generated
unit $R[F]$--module $N$, then $n_1,\ldots,n_t$ generate $N$ as a
$D_R$--module.
\end{corollary}

\begin{proof}
By Theorem \ref{thm.Main2} it is enough to check that the
$D_R$--submodule $M
\defeq D_R \cdot\langle n_1,\ldots,n_t \rangle$ satisfies $M \subseteq
F^*M$ and contains the $R[F]$--module generators $n_1,\ldots,n_t$
of $N$. The second statement is trivial and for the first one
observes that, by definition of root, one can write $n_i = \sum
r_j F(n_j)$ for some $r_j \in R$. Noting that $F(n_j) \in
F^*M$ we conclude $n_i \in F^*M$ for all $i$ as required.
\end{proof}

The above corollary is a generalization of Theorem \ref{thm.Main1}
in that $R_f$ is a finitely generated $R[F]$--module with root
generated by $n=\frac{1}{f}$.

\section{The case of a finitely generated algebra over a formal power
series ring}\label{formalsec}

The purpose of this section is to prove that $R_f$ is $D_{R|k}$-generated
by
$\frac{1}{f}$ in an important case that is not covered by our
previous results. Namely for $R$ a finitely generated algebra over
a power series ring $A=k[[x_1,\ldots,x_n]]$ over a field $k$ of
positive characteristic. The improvement is that we no longer assume
that $k$ is perfect.

Fixing the notation just introduced we further denote by
$k[[A^{p^s}]]=k[[x_1^{p^s},\dots, x_n^{p^s}]]$ the $k$-subalgebra
of $A$ consisting of all the power  series in $x_1^{p^s},\dots,
x_n^{p^s}$ with coefficients in $k$. By $k[[A^{p^s}]][R^{p^s}]$ we
denote the $k[[A^{p^s}]]$-subalgebra of $R$ generated by the
$p^s$-th powers of all the elements of $R$. The fact that $A$ is a
finite $k[[A^{p^s}]]$--module implies that $R$ is a finite
$k[[A^{p^s}]][R^{p^s}]$--module. Hence the ring of the
$k[[A^{p^s}]][R^{p^s}]$-linear differential operators of $R$ is
just $\text{End}_{k[[A^{p^s}]][R^{p^s}]}(R)$ due to formula (1) of
Section 2. Every $k[[A^{p^s}]][R^{p^t}]$-linear differential
operator of $R$ is automatically $k$-linear so $D_{R|k}\supseteq
V(R,k)$ where
\[V(R,k)=\bigcup_s\text{End}_{k[[A^{p^s}]][R^{p^s}]}(R)\]
As is pointed out in \cite[Ex. 5.3c]{Lyub}, we do not know whether
this  containment is always an equality (but it is if $R=A$).

Let $k^*=k^{1/p^{\infty}}$ be the perfect closure of $k$, let
$A^*=k^*[[x_1,\dots,x_n]]$ and let $R^*=A^*\otimes_AR$ where $A^*$
is regarded as an $A$-algebra via the natural inclusion
$k[[x_1,\dots,x_n]]\subseteq k^*[[x_1,\dots,x_n]]$. Since $R$ is a
finitely  generated $A$-algebra, $R^*$ is a finitely generated,
hence Noetherian, $A^*$-algebra.

\begin{theorem}\label{formal.thm}
With notation as above, let $R$ be a finitely generated
$A$-algebra such that $R^*$ is regular. Then $R_f$ is generated by
$\frac{1}{f}$ as a $D_{R|k}$--module.
\end{theorem}

\begin{proof} Since $V(R,k)$ is a subring of $D_{R|k}$, it is
enough to prove that $\frac{1}{f}$ generates $R_f$ as a
$V(R,k)$--module. According to \cite[p. 129]{Lyub},
$$D_{R^*}=R^*\otimes_RV(R,k)$$ and there is a functor
$$V(R,k)\text{-mod}\stackrel {N\mapsto R^*\otimes_R
N}{\longrightarrow} D_{R^*}\text{-mod}$$  where for each
$V(R,k)$--module $N$ the $D_{R^*}$--module structure on
$R^*\otimes_R N$ is defined as follows: if $\delta\in D_{R^*},
r\otimes n\in R^*\otimes_R N$ and $\delta (r)=\sum_i(r_i\otimes
v_i)$, where $r_i\in R^*, v_i\in V(R,k)$, then $\delta(r\otimes
n)=\sum_i(r_i\otimes v_i(n))$.

Since $A^*$ is flat over $A$ \cite[Thm. 22.3$(\beta)(1)(3')$]{Mat}
and local, it is faithfully flat over $A$, hence $R^*$ is
faithfully flat over $R$. Let $M\subset R_f$ be the
$V(R,k)$-submodule generated by $\frac{1}{f}$. Since
$R^*_f=R^*\otimes_RR_f$, we conclude that $R^*\otimes_RM$ is a
$D_{R^*}$-submodule of $R^*_f$ containing $\frac{1}{f}$ (we identify
$1\otimes f$ with $f$). By Theorem \ref{thm.Main1},
$R^*\otimes_RM=R^*_f$. Since $R^*_f$ is faithfully flat over $R$, we
conclude that $M=R_f$.
\end{proof}

\begin{theorem}\label{formal.gener}
With notation as above, let $R$ be a finitely generated
$A$-algebra such that $R^*$ is regular. If $n_1,\ldots,n_t$ are
generators of a root of a finitely generated unit $R[F]$--module
$N$, then $n_1,\ldots,n_t$ generate $N$ as a $D_{R|k}$--module.
\end{theorem}

\begin{proof} It is enough to prove that $n_1,\ldots,n_t$ generate
$N$ as a $V(R,k)$--module. If not, let $M\subseteq N$ be the
$V(R,k)$-submodule of $N$ generated by $n_1,\dots, n_t$. Since
$R^*$ is faithfully flat over $R$, we conclude that
$R^*\otimes_RM$ is a $D_{R^*}$-submodule of $R^*\otimes_RN$
different from $R^*\otimes_RN$. But this contradicts Corollary
\ref{coro.main} since $R^*\otimes_RN$ contains $1\otimes
n_1,\dots,1\otimes n_t$ and these elements generate a root of
$R^*\otimes_RN$. \end{proof}

The following special case of Theorems \ref{formal.thm} and
\ref{formal.gener} deserves to be stated separately.

\begin{corollary}
Let $R$ be a finitely generated algebra over a field $k$ of
characteristic $p>0$ such that $k^{1/p^{\infty}}\otimes_kR$ is regular.
Then

(a) $R_f$, for any $f\in R$, is generated by $\frac{1}{f}$ as a
$D_{R|k}$-module.

(b) More generally, if $n_1,\ldots,n_t$ are
generators of a root of a finitely generated unit $R[F]$--module
$N$, then $n_1,\ldots,n_t$ generate $N$ as a $D_{R|k}$--module.
\end{corollary}

\end{document}